\documentclass[aps,12pt,nofootinbib]{revtex4}

\usepackage{graphicx,graphics}
\usepackage{enumerate}
\usepackage{subfigure}
\usepackage{amsmath}
\usepackage{amsfonts}
\usepackage{amssymb}
\usepackage{amsthm}
\usepackage{psfrag}

\catcode`\@=11
\def\numberbysection{\@addtoreset{equation}{section}
\def\theequation{\arabic{section}.\arabic{equation}}}
\numberbysection

\def\Z{\mathbb{Z}}

\begin{document}

\title{Spiral Structures in the Rotor-Router Walk}

\author{Vl.V. Papoyan$^1$, V.S. Poghosyan$^2$ and V.B. Priezzhev$^1$}
\affiliation{
$^1$Bogoliubov Laboratory of Theoretical Physics,\\ Joint Institute for Nuclear Research, 141980 Dubna, Russia\\
$^2$Institute for Informatics and Automation Problems\\ NAS of Armenia, 0014 Yerevan, Armenia
}

\begin{abstract}
We study the rotor-router walk on the infinite square lattice with the outgoing edges at each lattice site ordered clockwise.
In the previous paper [J.Phys.A: Math. Theor. 48, 285203 (2015)], we have considered the loops created by rotors and labeled sites where the loops become closed. The sequence of labels in the rotor-router walk was conjectured to form a spiral structure obeying asymptotically an Archimedean property.
In the present paper, we select a subset of labels called ``nodes'' and consider spirals formed by nodes.
The new spirals are directly related to tree-like structures
which represent the evolution of the cluster of vertices visited by the walk.
We show that the average number of visits to the origin $\left<n_0(t)\right>$ by the moment $t\gg 1$ is $\left<n_0(t)\right> = 4 \left<n(t)\right> + O(1)$ where $\left<n(t)\right>$ is the average number of rotations of the spiral.
\end{abstract}

\maketitle

\noindent \emph{Keywords}: rotor-router walk, Archimedean spiral, sub-diffusion.

\section{Introduction}

The rotor mechanism, firstly proposed in the theory of self-organized criticality \cite{BTW,Dhar} under name ``Eulerian walk''
\cite{PDDK}, was rediscovered independently as a tool for a derandomization of the random walk \cite{CS,HP}.
The subsequent studies were concerned with collective properties of the medium ``organized'' by the walk and with statistical properties of the walk itself \cite{LP05,LP07,LP08,AngelHol,HussSava,FGLP,FLP}.

\begin{figure}[!ht]
\includegraphics[width=160mm]{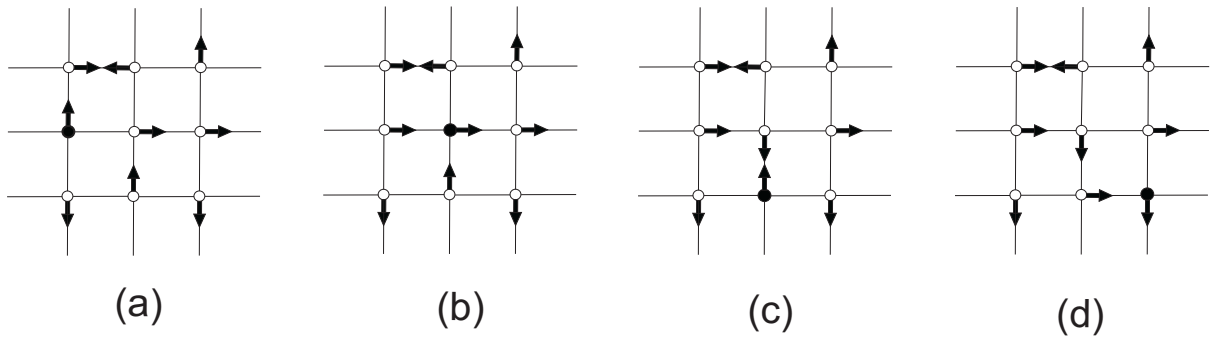}
\caption{Circles denote the lattice sites. (a) The chip is originally in the filled circle where the arrow is directed ``up''.
(b) The chip rotates the arrow clockwise and moves right. (c) The next clockwise rotation sends the chip down. (d) The last position
of the chip is the right lower corner.}
\label{steps}
\end{figure}

The dynamics of the rotor-router walk can be described as follows.
Consider a square lattice with arrows attached to the lattice sites.
Arrows attached to the lattice sites are directed toward one of their neighbors on the lattice.
A particle called usually {\it chip}, performs a walk jumping from a site to a neighboring site.
Arriving to a given site, the chip changes direction of the arrow at that site in a prescribed
order and moves toward the neighbor pointed by new position of the arrow.
Thus, given an initial orientation of arrows on the whole lattice, the rotor-router walk is deterministic.
The walk started from uniformly distributed random initial configurations can be called uniform rotor walk.
Three steps of the rotor walk on the square lattice are shown in  Fig.\ref{steps}.

If the lattice is finite, the walk starting from an arbitrary site settles into an Eulerian circuit where each edge of the lattice is visited exactly once in each direction  \cite{PDDK,HLMPPW}.
When the walker is in the Eulerian circuit, configurations of rotors $\rho$ associated to each site are recurrent.
A graphic representation of the recurrent configuration is {\it unicycle} which is a specific state
where the arrows form a spanning set of directed edges containing a unique directed cycle which the chip belongs to \cite{HLMPPW}. If the position of the chip on the cycle is $a$, we denote the unicycle as $(\rho,a)$.

Along with the unicycle, we can define the {\it multicycle} \cite{PPP} as a graph containing exactly $k$ cycles together with $k$ chips at vertices $a_0,a_1,\dots,a_{k-1}$ belonging to the cycles. For multicycles, we use the notation $(\rho,a_0,a_1,\dots,a_{k-1})$.

For the infinite lattice, both questions on trajectories of the walker and on the configurations of arrows become more complicated.
A basic problem here is to find the range of the rotor walk, i.e. the number of distinct sites visited in $t$ steps and,
given the lattice symmetry and the rotor mechanism, to find a shape of the area visited by the walker.
One conjecture and one theorem shed light on this problem.
Kapri and Dhar \cite{Kapri} conjectured that the set of sites visited by the clockwise uniform rotor walk in $t$ steps on the infinite square lattice is
asymptotically a disk of average radius $ct^{1/3}$ where $c$ is a constant.
Florescu, Levine and Peres \cite{FLP} proved that for an infinite $d$-dimensional lattice, regardless of a rotor mechanism or an initial rotor configuration,
the rotor walk in $t$ steps visits at least on the order of $t^{d/(d+1)}$ distinct sites.

Monte Carlo simulations in \cite{Kapri} showed that the average number of visits of a site inside the disk is a linear decreasing function of its distance from the origin.
The authors of \cite{Kapri} give the following explanation of this characteristic behavior.
After a moment when two sites at different distances from the origin are visited by the rotor walk,
both sites are visited equally often because of the local Euler-like organization of arrows.
Then, the difference between the numbers of visits of these sites remains bounded for an arbitrary number of subsequent steps.
The key point in this explanation is the local Eulerian  organization which is proven rigorously only for finite graphs as a part of the total organization.
For the infinite lattice, any bounded domain tends to the entirely organized state only asymptotically being repeatedly visited by the rotor walk.
A question, however, is in the periodicity of returns.
The mean number of returns and the mean-square displacement should be in a definite proportion to provide the sub-diffusive behavior of the rotor walk.
So, it is desirable to find in the system of rotors some structure which provides sufficiently often returns of the walker to the origin and, as a consequence, to any previously visited site.
The construction of such a structure is the main goal of the present paper.

In the recent work \cite{PPP}, we have considered the motion of the clockwise rotor-router walk inside closed contours emerged in random rotor configurations on the infinite square lattice. We proved a property called the {\it weak reversibility}: even though the configuration of rotors inside the contour is random, the rotor-router walk inside the contour demonstrates some regularity, namely, the chip entering the clockwise contour $C$ in a vertex $v \in C$ leaves the contour at the same vertex $v$, and then the clockwise orientation of rotors on $C$ becomes anti-clockwise.

We referred to the sites where rotors complete clockwise contours as {\it labels}, and noticed that the sequence of labels forms a spiral structure. After averaging over initial random configurations of rotors, the sequence approaches asymptotically the {\it Archimedean} spiral. However, the spiral structure as such does not explain the obligatory periodic visits of the origin by the rotor walk. In Section III, we consider particular labels called {\it nodes}. The set of nodes being a subset of that of labels has also the spiral structure. The difference between labels and nodes lies in the disposition of contours corresponding to them. In the case of labels, a contour completed at given site is not necessarily adjacent to the contour associated with the previous label. In case of nodes, each new contour associated with a node either has common sites with that corresponding to the previous node, or contains this contour inside.

In Section IV, we analyze the structure of contours associated with nodes.
According to the week reversibility, each contour after visiting its interior becomes anti-clockwise and left by the chip.
At the moment of exit from the node located at vertex $v$, there is a directed path formed by arrows and connecting one of
neighboring sites of $v$ with the previous node.
A collection of paths obtained at the moment $t$ is a tree rooted at the current location of the chip.
The tree structure together with spiral-like motion of the chip provides the obligatory visits to the origin for each turn of the spiral.
Depending on the location of the origin with respect to the first clockwise contour, the chip returns to the vicinity of origin by different ways.
Once the spiral structure is formed, the number of visits to the origin is 4 for each rotation around the origin.
Then, the total average number of visits $\left<n_0(t)\right>$ for $\left<n(t)\right>$ rotations performed by the chip starting from the uniform random initial configurations of arrows is $\left<n_0(t)\right> = 4 \left<n(t)\right> + O(1)$, when $t \gg 1$.
Since the Archimedean spiral has a constant interval between coils, we obtain the linear dependence between the radius of the spiral and the number of returns to the origin.

In the separate section V, we analyze the convergence of the set of labels (nodes) to the Archimedean spiral.
As it was noticed in \cite{PPP}, this convergence is extremely slow.
The existence of the limit for the averaged ratio of radius $r$ to angle $\theta$ which is a constant $b$ for the purely Archimedean case, is not proven yet. A prospective value of $b$ can be obtained from the scaling law for the average number of
visits conjectured by Kapri and Dhar \cite{Kapri}. Extensive simulations show that the deviation from the constant $b$ remains considerable for very large number of nodes $7\cdot 10^5$ and huge number of steps $10^{10}$. We discuss a possible reason for this deviation.

\section{Contours in rotor state and chip motion}

We consider the infinite square lattice, and fix the clockwise rotor mechanism at each site.
In the initial rotor state, the arrows at each lattice site are directed randomly to one of four directions with equal probabilities, and the chip is in the origin.
At each step of discrete time, the chip arriving at a site rotates the arrow at that site 90 degrees clockwise, and moves to the neighboring site pointed by the new position of arrows.

The motion of the chip is determined by the current rotor state. Given a rotor state, we say that a group of arrows forms a directed path if the arrows are attached to sites $v_1,v_2, \dots, v_{n}$ such that $v_i$ and $v_{i+1}$ are neighbors, and the arrow at $v_i$ is directed toward $v_{i+1}$ for $i=1,\dots,n-1$.
The directed path of arrows becomes a cycle if  $v_1 = v_{n}$.
A shortest possible cycle consists of two adjacent sites $v_1$, $v_2$, which are connected by a pair of edges from $v_1$ to $v_2$ and back. We call such cycles {\it dimers} by analogy with lattice dimers covering two neighboring sites. A cycle formed by more than two edges is called {\it contour}.

The configuration of arrows inside a contour is either free of cycles or contains a number of cycles.
In the first case, the arrows inside the contour form a spanning forest rooted at the contour.
In the second case, the arrows form a spanning forest where each tree is rooted either at the external contour, or at one of the internal cycles. Correspondingly, two theorems describe the behavior of the chip approaching contours.

{\it Theorem 1} on reversibility (\cite{PPS} and \cite{HLMPPW}, Corollary 4.11).
Let $G$ be a planar graph containing a unicycle $(\rho,a)$ with the contour $C$ oriented clockwise and $a\in C$.
After the rotor-router walk makes some number of steps, each rotor internal to $C$ has performed a full rotation,
each rotor external to $C$ has not moved, each rotor on $C$ has performed a partial rotation so that $C$ is now oriented anti-clockwise and the chip has returned to $a$.

To describe the motion of the chip in the second case, we consider a more general situation, where several chips are involved
into the evolution of arrows inside the contour.

{\it Theorem 2}, \cite{PPP}. Let $G$ be a connected bidirected planar graph and $(\rho,a_0,a_1,\dots,a_{k-1})$ be a multicycle with the external contour $C_0$ oriented clockwise together with $k-1$ internal cycles $C_1,\dots,C_{k-1}$ oriented anti-clockwise. The rotor-router operation is sequentially applied to the chip at $a_0 \in C_0$ until the moment $T_0$ when the chip returns to $a_0$, and the rotor at $a_0$ is made oriented anticlockwise.
Then, the same is applied to chips at $a_1,\dots,a_{k-1}$ until the moments $T_i$ when chips starting from $a_i \in C_i$ return to $a_i$ and the rotors at $a_i$ are made oriented clockwise.
Then, all rotors on $C_0$ are becoming oriented anticlockwise, all rotors on $C_1,\dots,C_{k-1}$ become oriented clockwise, and all vertices internal to $C_0$ and external to $C_1,\dots,C_{k-1}$ perform a full rotation.

The description of motion of the single chip inside a contour is given by a reduced version of Theorem 2 :

{\it Theorem 3} on weak reversibility. Let $G$ be a planar graph containing the external contour $C_0$ oriented clockwise and some number of internal cycles inside $C_0$. The rotor-router walk starting at site $v$, $v \in C_0$ moves until the moment when the chip returns to $v$. As a result, all rotors on $C_0$ become oriented anti-clockwise, and the chip leaves contour $C_0$ at the next time step. The rotors internal to $C$ perform either a full rotation or a partial rotation or do not move at all.

The Theorem 2 allows us to specify internal sites of contour $C_0$ which perform the full rotation.

{\it Corollary.} Let $F$ be a set of sites inside the clockwise contour $C_0$ which belong to the forest rooted at $C_0$. Assume that there is the site $v_{in} \in F$ and  $v_{in} \ni C_0$. The chip starting at site $v \in C_0$ moves until the moment when the chip returns to $v$. Then, the rotor at $v_{in}$ performs the full rotation.

{\it Proof.} Let $T_{tot} = \max_{i\leq k-1} (T_i)$ be the total number of steps in the process described in theorem 2. The process can be divided into two stages, I and II: stage I for steps $0\leq T \leq T_0$ and stage II for steps $T_0 < T\leq T_{tot}$. Consider the site $v_{in}\in F$ belonging to the forest $F$ rooted at $C_0$. According to theorem 2, the rotor at $v_{in}$ performs the full rotation to the moment $T_{tot}$. The number of rotations is $\deg (v_{in})$. Assume that a part of rotations is performed during stage I and the rest of them during stage II. Since stage II follows stage I, a sequence of arrows resulting after $T_{tot}$ steps is directed from $v_{in}$ to one of sites $a_i, 1\leq i \leq k-1$. But $v_{in}$ belongs to a tree rooted at $C_0$ by the condition. Therefore,all $\deg (v_{in})$ rotations are performed during stage I, which coincides with the process described in Theorem 3.

Below, we apply the Theorems 1,3 and Corollary to investigate the uniform rotor walk on the infinite square lattice.

\section{Labels, nodes and spirals}

\begin{figure}[!ht]
\includegraphics[width=160mm]{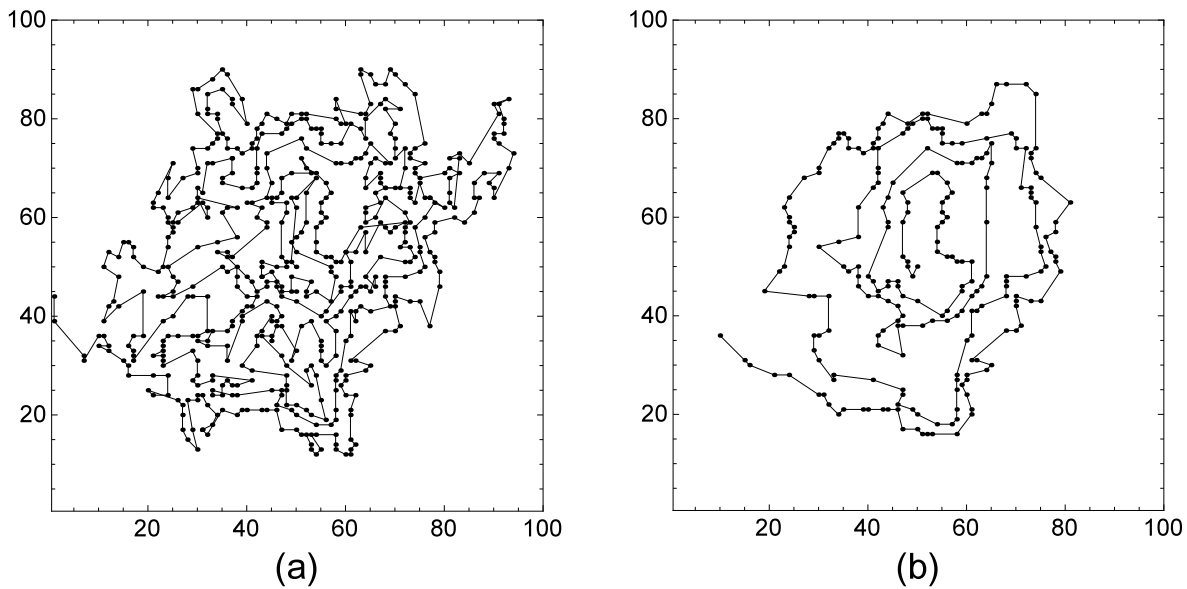}
\caption{(a) Spiral of labels; (b) spiral of nodes.}
\label{fig2}
\end{figure}

The rotor walk during the time evolution creates contours of arrows $C_1,C_2,\dots$ sequentially. In \cite{PPP}, we considered the set of sites $v_1,v_2,\dots$ where each contour becomes closed at time steps $t_1,t_2,\dots$, and called these sites labels.
According to theorems 1,3, the rotor walk leaves each contour $C_i$ at step $t_i^{'} > t_i$ at the same site $v_i,i=1,2\dots$. We skip details of evolution in the time intervals between $t_i$ and $t_i^{'}$ and ignore possible new contours appearing inside $C_i$ during these intervals. The only fact of the evolution between $t_i$ and $t_i^{'}$ we take into account is reversing the
clockwise orientation of the contour $C_i$. It was found in \cite{PPP} that the labels $v_1,v_2,\dots$ are not simply situated in the cluster of visited sites but form a spiral structure. An example of the spiral of labels is shown in Fig.\ref{fig2}(a).

Whereas each particular spiral has an irregular form, their average over uniform random initial rotor states tends to the spiral obeying the Archimedean property
\begin{equation}
r=a+b\,\theta,
\label{Archimed}
\end{equation}
in planar coordinates $r,\theta$, with constant $a$ and $b$. A convergence of spirals to the Archimedean law (\ref{Archimed}) for a large number of steps is discussed in Section V.

Despite the surprising property, the spiral of labels does not say anything about periodicity of the chip returns to the vicinity of the origin needed for organization of the cluster of visited sites. To answer this question, we consider the sequence of contours $C_1,C_2,\dots$ corresponding to labels in more detail.

Let $C_i$ and $C_{i+1}$ be two successive contours in the sequence $C_1,C_2,\dots$. There are three possibilities for the disposition of $C_{i+1}$ with respect to $C_{i}$: a) the set of sites $\{v\}_{C_{i+1}}$ where arrows of $C_{i+1}$ are attached
has no common sites with $\{v\}_{C_{i}}$ and contour $C_{i+1}$ is outside $C_{i}$; b)the set $\{v\}_{C_{i+1}}$ has no common sites with $\{v\}_{C_{i}}$ and contour $C_{i}$ is inside $C_{i+1}$; c) the set $\{v\}_{C_{i+1}}$ has at least one common site with $\{v\}_{C_{i}}$. To provide for the condition b), the contour $C_{i+1}$ should contain inside at the moment $t_{i+1}$ all sites visited at moments $t\leq t_{i}$. Otherwise, there are lattice sites outside $C_{i+1}$ which do not connected with $C_i$ at the moment $t_i$ by any path of arrows, what is impossible for a single walk. When the cluster of visited sites grows, the probability of a contour enveloping the cluster of  previously visited sites dramatically decreases and we can exclude the case b) from the consideration. Then,we select from the set of labels $v_1,v_2,\dots$ a subset of labels $v_{i_1},v_{i_2},\dots$ whose contours obey criterion c) and assume that $v_{i_1}$ coincides with $v_1$. We call the selected labels {\it nodes}. Fig.\ref{fig2}(b) shows the spiral consisting of nodes selected from the labels of Fig.\ref{fig2}(a).

Let $\bar{v}_i, i=1,2,\dots$ be the sequence of nodes generated by the rotor walk. By the construction,
at the moment of exit from the node located at vertex $\bar{v}_{i+1}$, there is a directed path formed by arrows and connecting one of neighboring sites of $\bar{v}_{i+1}$ with the previous node $\bar{v}_{i}$. Thus, the system of directed paths and remaining parts of contours associated with nodes forms a connected network. In what follows, we will see that the obtained network is a tree constructed from topologically uniform elements.

\section{Tree structure and number of visits the origin}

\begin{figure}[!ht]
\includegraphics[width=50mm]{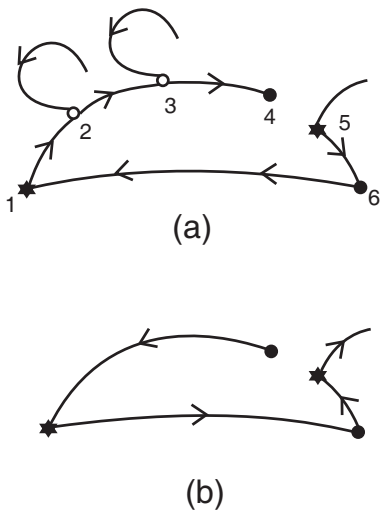}
\caption{(a) The rotor sequence traced by chip moving from node 1 to node 5 (see text). (b)The building block of the tree. }
\label{LabStar}
\end{figure}

The chip moving from a fixed site $v_a$ to site $v_b$ traces a path of arrows directed from $v_a$ to $v_b$. Fig.\ref{LabStar}(a) illustrates the path traced by the chip moving between two successive nodes $\bar{v_{i}}$ and $\bar{v}_{i+1} $ selected from the sequence of nodes $\bar{v}_i, i=1,2,\dots$. The first node is at the site $\bar{v_{i}}$ marked by 1. The directed path between sites 6 and 1 is a part of the contour $C(\bar{v_{i}})$ corresponding to the first node. The chip creates a clockwise contour at site 2, reverses it to the anti-clockwise one and leaves it at the same site 2. Assume, that this contour has no common sites with line $(6,1)$. If it also has no intersections with either of previous contours, the site 2 is a label but not a node. The situation is repeated at site 3 and the chip continues motion creating an arbitrary number of labels until reaching site 5. In general, site 5 does not belong to line $(6,1)$ because of a possible sequence of arrows directed from site 5 to site 6 existing before the chip could reach it. As a result, the clockwise contour appears $C(\bar{v}_{i+1})=(1,2,3,4,5,6,1)$ at site 5 which is a node $\bar{v}_{i+1}$ according to criterium c). After reversing $C(\bar{v}_{i+1})$ to the anti-clockwise contour $\bar{C}(\bar{v}_{i+1})=(5,4,3,2,1,6,5)$, the chip leaves $\bar{C}(\bar{v}_{i+1})$ at site 5. The parts of reversed contours corresponding to the labels 2 and 3 and possible others are just branches attached to the path from 4 to 5 which do not affect connectivity of this path and therefore can be ignored. The resulting path shown in
Fig.\ref{LabStar}(b) is a building block of the tree we are going to construct.

\begin{figure}[!ht]
\includegraphics[width=60mm]{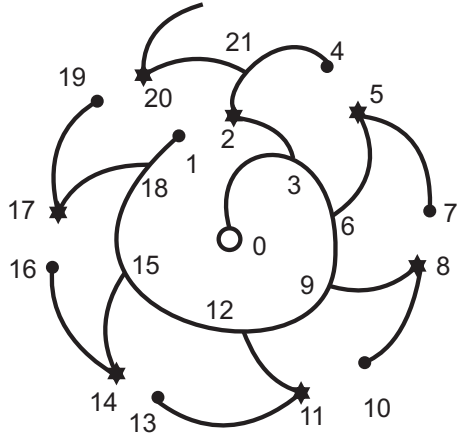}
\caption{The tree formed after creation of nodes 2,5,8,11,14,17,20 (see text).}
\label{Tree}
\end{figure}

The tree is constructed by consecutive adding the building blocks, one by one for each new node. In Fig.\ref{Tree} we show
schematically how the tree grows if one neglects a possible difference in sizes of the blocks and takes into account their topological structure only.
Consider the first clockwise contour containing the origin. In  Fig.\ref{Tree}, this contour is formed by the path starting at site $0$ and reaching site 2. Due to the possible directed sequence from site 2 to 3, the contour is closed at site 2 which is the first node. At site 2, the chip leaves the contour after reversing its direction and reaches site 5 which is connected with site 6 by the directed sequence of arrows, if exists (otherwise site 5 coincides with 6). Then the clockwise contour 2,4,5,6,3,2 appears with the next node at site 5. The interior of this contour covers the sector of space bounded by two radial branches 3,2,4 and 3,6,5. Reversing this contour to the anti-clockwise one, the chip continues motion from node 5 to node 8 and so on. By the construction of the building blocks, the resulting graph $T$ is a tree. At each moment of time, all directed paths of arrows
constituting the tree are oriented towards the current position of the chip.
Consider the chip at the last node situated at site 20. Its further path can trace a directed path of arrows to any reachable site of the tree $T$. Assume this site is on the interval between sites 10 and 8. Then the next clockwise contour will contain the sites 8,9,6,3,2,21,20 and cover the sector bounded by radial branches 3,2,21,20 and 3,6,9,8. The continuation of the spiral up to full rotation creates subsequent contours covering one by one {\it all} sectors between branches of the tree, including the sector bounded by line $20,21,2,\dots,19$ and containing the origin.

If the initial rotor configuration contains one, two or three sequences of arrows flowing into the origin, the structure of the tree $T$ admits two,three or four sectors having the origin on the boundary between sectors.

The described scenario being empirical is nevertheless typical for any random initial configuration of rotors. The obtained tree is stable after each rotation of the spiral. Indeed, by theorem 2, the chip visiting a contour changes its orientation but not its form. Branches of the tree situated inside contours consist of sites which perform the full rotation according to the Corollary in section II and therefore remain stable as well. The arrows inside the contours not belonging to these branches cannot create new cycles due to the rotor mechanism and can only add new branches to the existing tree.

The main conclusion we can draw from the existence of the tree structure and the spiral ordering of nodes is that every turn of the spiral generates necessarily either a contour containing the origin inside it or generates contours containing the origin on their boundaries. Now we are ready to answer the question about the number of visits to the origin for each rotation.

\begin{figure}[!ht]
\includegraphics[width=70mm]{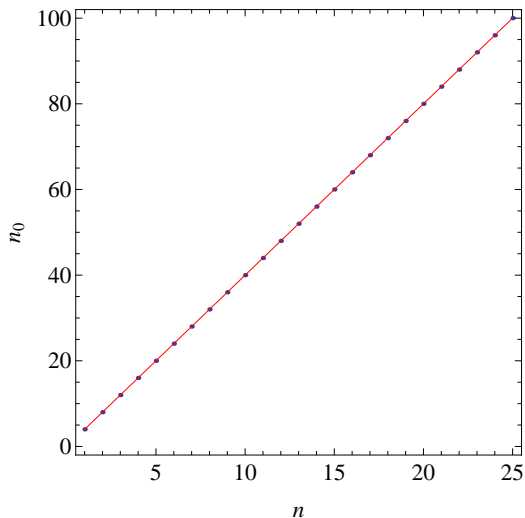}
\caption{Number of visits $n_0$ vs number of rotations $n$.}
\label{s1}
\end{figure}

Consider the first clockwise contour $C_{first}$ containing the origin inside. By the construction of the tree $T$, the origin $0$ belongs to the forest rooted at sites of $C_{first}$. It follows from the Corollary, that the arrow at the origin performs the full rotation when $C_{first}$ is reversed. As the tree $T$ grows, the subsequent contours containing the origin become larger and the chip visiting them not necessarily visits all sites of their interiors. However, the forest arisen in the first contour remains rooted at all subsequent contours containing the origin. Then by the Corollary, the rotor at the origin performs the full rotation each time when such a contour appears. Since it happens for each turn of the spiral, the number of visits of the origin $n_0$ depends of the number of turns $n$ as
\begin{equation}
n_0= 4 n+ O(1)
\label{number}
\end{equation}
where $O(1)$ is due to possible visits to the origin before the moment when the first loop of the spiral is formed.
The result of simulation for a single spiral is shown in Fig.\ref{s1}.

If the initial configuration contains one or more sequences of arrows directed to $0$,the origin can belong to the boundary of adjacent sectors of the tree. Since all sectors are covered by contours appearing during one rotation of the spiral, the total number of visits to the origin remains 4. For instance, let $C_1$ and $C_2$ be two adjacent contours and the origin 0 is at the boundary between them. Reversing contour $C_1$ at some stage of rotation, the chip rotates the arrow at 0 by angle $\pi/4, \pi/2$ or $3\pi/4$ depending on the form of the boundary. Contour $C_2$ is also reversed during the same turn with the arrow rotation by the additional angles $3\pi/4, \pi/2$ or $\pi/4$ giving the full rotation per one turn.

\section{Convergence to the Archimedean spiral}

The emergence of a spiral structure in the clockwise rotor-router walk comes from a rather simple reason. Indeed,consider a label $v_k$ situated near the boundary of the cluster of visited sites. A preferable position for the next label $v_{k+1}$ is on the right side from the branch $0 \rightarrow v_k$ to provide the clockwise orientation of the contour $C_{k+1}$ if it has common edges with $0 \rightarrow v_k$. Then, the preferable direction of successive positions of labels $v_k,v_{k+1},v_{k+2},\dots$ is clockwise with respect to the origin of the cluster. Since the size of the cluster grows with time, the positions $v_k,v_{k+1},v_{k+2},\dots$ form a spiral-like structure.

The set of nodes, being a subset of labels has the spiral form as well. Moreover, the condition for contours corresponding to two successive nodes to be adjacent imposes an additional restriction on the positions of nodes. The comparison of typical spirals in Fig.\ref{fig2}(a) and (b) shows that this restriction makes the spiral of nodes more regular than that of labels.

Introducing the polar coordinates, we denote by $r(k)$ the distance from the origin
and by $\theta(k)$ the winding angle of $k$-th node. Then, we say that the spiral of nodes is asymptotically Archimedean in average if
\begin{equation}
\left\langle \frac{r(k)}{\theta(k)} \right\rangle \rightarrow b \quad \text{for} \quad k\rightarrow \infty,
\label{definition}
\end{equation}
where the average is taken over the uniformly distributed  states of the spiral.

A numerical verification of the Archimedean property for labels in \cite{PPP}  showed a very slow converge to the asymptotic law (\ref{definition}). Here, we use the more pronounced spiral structure of nodes to determine the asymptotic behavior of $r(k)/\theta(k)$ with greater accuracy.

First, we compare the conjectured Archimedean property (\ref{definition}) with the conjecture by Kapri and Dhar \cite{Kapri},
who supposed that the average number of visits $n_N(x)$ to the site separated from the origin by distance $x$ for $N$
steps satisfies the scaling form
\begin{equation}
n_N(x)=a N^{1/3}F\left(\frac{x}{cN^{1/3}}\right),
\label{scaling}
\end{equation}
where $F(y)$ is the scaling function $F(y)=1-y$, $0\leq y < 1$.
The data obtained in \cite{Kapri} for $N=10^6, 10^7$ and $10^8$ collapse with $a=0.5$ and $c=1.38$.
Then, the limiting radius of the circle of visited sites
depends on the number of visits to the origin  $n_0$ as
\begin{equation}
r_{circle}=\frac{c}{a}n_0 \simeq  2.76 n_0
\label{Dharrad}
\end{equation}

The spiral of nodes, by the construction, is contained inside the circle of visited sites. According to the Archimedean law
Eq.(\ref{Archimed}) with $\theta =2 \pi n$, for the large number of rotations $n$, the spiral curve tends to the circle of
radius $r_{spiral}= 2\pi b n$. Then, due to Eq.(\ref{number}), we obtain again the linear dependence of the radius on the number of visits to the origin
\begin{equation}
r_{spiral}=\frac{b\pi n_0}{2}
\label{PPPrad}
\end{equation}
Thus, the conjectured Archimedean property is consistent with the linear scaling law Eq.(\ref{scaling}). In the scaling limit, we can expect that $r_{spiral}=r_{circle}$. Then, the constant $b$ in Eq.(\ref{definition}) is related to the constants of the scaling law as $b=2c/(a\pi)$. Taking $a=0.5$ and $c=1.38$, we obtain $b=1.757$.

\begin{figure}[!ht]
\includegraphics[width=90mm]{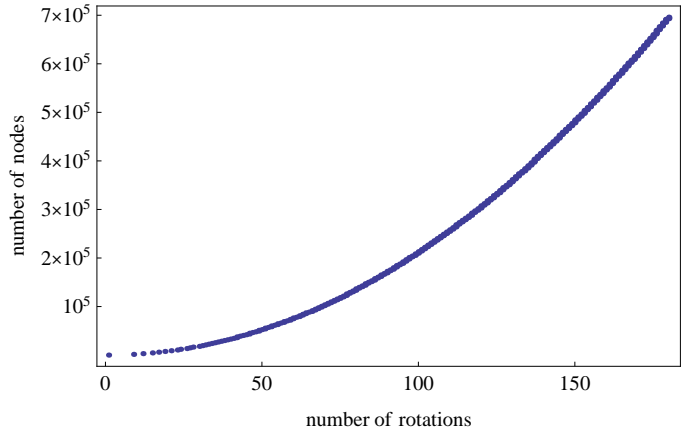}
\caption{Number of nodes vs number of rotations.}
\label{fig5}
\end{figure}

To study the asymptotical behavior of the spiral, we extended our simulations to number of time steps $T \sim 10^{10}$ and number of samples $200000$. The numbers of nodes and spiral rotations we use in our analysis are shown in Fig.\ref{fig5}.

\begin{figure}[!ht]
\includegraphics[width=110mm]{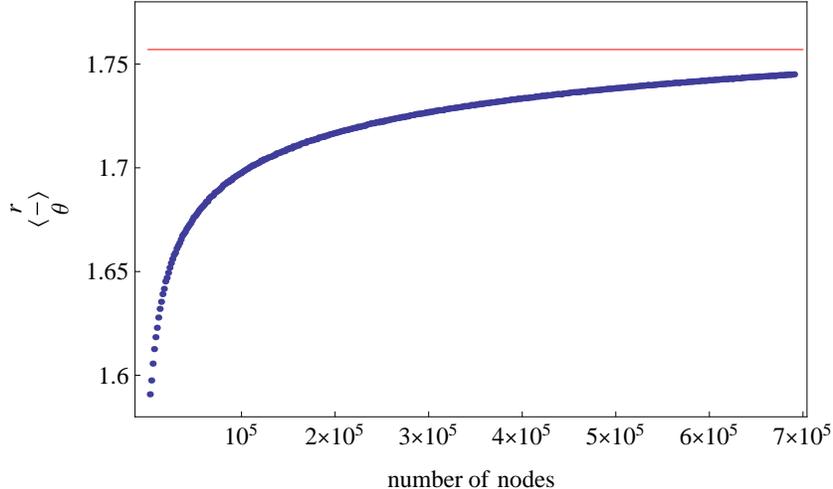}
\caption{Average ratio $r/\theta$ vs number of node $k$.}
\label{fig4}
\end{figure}

The result of averaging of the ratio $r/\theta$ over uniformly distributed initial conditions as a function of the node number
is shown in Fig.\ref{fig4}. The horizontal line in Fig.\ref{fig4} corresponds to $b=1.757$ obtained from the scaling conjecture
\cite{Kapri}. We can see that the deviation from the constant $b$ remains considerable up to very large node numbers.
Despite the apparent decrease of the slope of the curve in Fig.\ref{fig4}, we cannot guarantee the coincidence of the limiting value of $\langle r/\theta \rangle$ with the scaling value $b=1.757$ and even the existence of the limit as such.

The Monte-Carlo simulations show that the ratio of radius to angle grows with $k$ not faster than $\log(k)^{\alpha}$ with $\alpha=0.16$.
The logarithmic deviation obtained for finite $k$ suggests
to try an asymptotic expansion in powers of $1/\log(k)$.
The obtained data yield the following lower bound:
\begin{equation}
<\frac{r(k)}{\theta(k)}> \simeq 2.08 - \frac{5.11}{\log k} + \frac{8.08}{\log^2 k}.
\end{equation}

A reason for so slow convergence is a geometrical non-equivalence of averaging spirals. Indeed, if the random spirals differ one
from another only by the distances between coils, the averaging over large number of samples would give a well defined mean distance even for a relatively broad distance distribution. Instead, we observe some number of meanders in the spiral structure at different regions of the spiral. This leads to an effective broadening of the intervals between coils, but since the meanders are rare events, the determination of parameters of the averaged spiral needs an enormously large statistics.

\section*{Acknowledgments}
VBP thanks the RFBR for support by grant 16-02-00252.
VSP thanks the JINR program ``Ter-Antonyan - Smorodinsky''.

\end{document}